\documentclass[11pt,a4paper]{article}
\usepackage{mathrsfs}
\usepackage{mathrsfs}
\usepackage{mathrsfs}
\usepackage{amsfonts}
\usepackage{amsfonts}
\usepackage{amsfonts}
\usepackage{mathrsfs}
\usepackage{amssymb}
\usepackage{amsmath}
\newtheorem{theorem}{Theorem}[section]

\numberwithin{equation}{section}

\begin{document}
\parindent 9mm
\title{Stabilization of parallel-flow heat exchangers with arbitrary delayed boundary feedback
\thanks{This work was supported by the Fundamental Research Funds for the Central Universities (Grant No. xjj2017177), the National Natural Science Foundation of China (grant no. 11301412), Natural Science Foundation of Shaanxi Province (grant no. 2014JQ1017),  Research Fund for the Doctoral Program of Higher Education of China (grant no. 20130201120053), Project funded by China Postdoctoral Science Foundation (grant nos. 2014M550482 and 2015T81011).}
}
\author{Zhan-Dong Mei \thanks{School of Mathematics and Statistics, Xi'an Jiaotong
University, Xi'an 710049, China; Email: zhdmei@mail.xjtu.edu.cn }}


\date{}
\maketitle \thispagestyle{empty}
\begin{abstract}
%
%
This paper is concerned with the stabilization of
parallel-flow heat exchangers equation with boundary control
and arbitrarily delayed observation. The observer and predictor-based method,
developed by [B.Z. Guo, K.Y. Yang, Automatica, 2009, 45(6): 1468-1475], is used to
stabilize equation. The exponentially stabilization of such system by estimated state
feedback control law is proved.
\vspace{0.5cm} 

%
%

\end{abstract}


\section{Introduction and main results}

The parallel-flow heat exchangers, which facilitates transfer of heat of two fowing parallel to each other in a coaxial tube, is
described as follows:
\begin{align}\label{delay}
    \left\{
      \begin{array}{lll}
        \frac{\partial}{\partial t}\theta_1(t,x)=-\frac{\partial}{\partial x}\theta_1(t,x)+h_1(\theta_2(t,x)-\theta_1(t,x)), 0<x<l, t>0, & \hbox{ } \\
        \frac{\partial}{\partial t}\theta_2(t,x)=-\frac{\partial}{\partial x}\theta_2(t,x)+h_2(\theta_1(t,x)-\theta_2(t,x)), 0<x<l, t>0, & \hbox{ } \\
        \theta_1(t,0)=u_1(t),\theta_2(t,0)=u_2(t), t\geq 0,& \hbox{ } \\
        y_1(t)=\theta_2(t-\tau,l), y_2(t)=\theta_1(t-\tau,l),  t\geq\tau, & \hbox{ } \\
        \theta_1(0,x)=\theta_{10}(x), \theta_2(0,x)=\theta_{20}(x),
      \end{array}
    \right.
\end{align}
where $l$ is length of two tubes, $\theta_1(t,x)$, $\theta_2(t,x)\in R$ are the temperature variations at time $t$ and at the point $x\in[0,l]$,
 $h_1, h_2>0$ denote the heat exchange rate, $u_1(t)$ and $u_2(t)$ are the boundary control, $\tau>0$ is a (known) constant time delay, and $y_1(t)$ and $y_2(t)$ are the observation which suffers from the time delay $\tau>0$.

The author in \cite{Sano2016} proved that system (\ref{delay}) with $u_1(t)=0$ and $u_2(t)=-ky_2$ is exponentially stable whenever $k^2<\frac{h_2}{h_1}$ and $h_1l<\tau<\frac{h_2l}{k^2}$. However, the exponentially stability of the cases $\tau<h_1l$ and $\tau>\frac{h_2l}{k^2}$ are unknown even when $k_1=0$, see Conclusion of \cite{Sano2016}. {\it Motivated by this}, in this paper, we shall use the scheme observer and predictor-based, developed by Guo and Yang \cite{Guo2009}, to stabilize the equation with arbitrarily delayed observation.

Let $A=A_0+A_1$ with
$A_0=\left(
\begin{array}{cc}
 -\frac{\partial}{\partial x} & 0\\
  0 & -\frac{\partial}{\partial x}\\
  \end{array}
  \right)$, $A_1=\left(
            \begin{array}{cc}
              -h_1 & h_1 \\
              h_2 & -h_2 \\
            \end{array}
          \right),$ $D(A)=\bigg\{\left(
\begin{array}{cc}
 f\\
  g\\
  \end{array}
  \right)\in H^1(0,l)\times H^1(0,l): f(0)=g(0)=0\bigg\}.$ By \cite{Chen2010}, $A$ generates a $C_0$-semigroup denoted by $\{e^{At}\}_{t\geq 0}$.
Denote by $B$ and $C$ the control operator and observation operator of the delay free system ($\tau=0$) corresponding to (\ref{delay}), respectively.
The corresponding state space is $L^2([0,1])\times L^2([0,1])$, input space and observation space are $R^2$. We note by $\|z\|$ the norm of $z$ on the associated Hilbert space. By \cite{Chen2010}, $B$ is admissible for $A$, that is,
for zero initial value current state is depended continuously on the input with $L^2$ norm.
Through simple computation, it follows that system
\begin{align*}
    \left\{
      \begin{array}{ll}
        \frac{\partial}{\partial t}\theta_1(t,x)=-\frac{\partial}{\partial x}\theta_1(t,x), 0<x<1, t>0 & \hbox{ } \\
        \frac{\partial}{\partial t}\theta_2(t,x)=-\frac{\partial}{\partial x}\theta_2(t,x), 0<x<1, t>0 & \hbox{ } \\
        \theta_1(t,0)=0,\theta_2(t,0)=0, t\geq 0,& \hbox{ } \\
        \theta_1(0,x)=\theta_{10}(x), \theta_2(0,x)=\theta_{20}(x),    0<x<1,& \hbox{ }
      \end{array}
    \right.
\end{align*}
implies
$\int_0^l|\theta_1(t,l)|^2dt=\int_0^l|\theta_{10}(l-t)|^2dt\leq\|\theta_{10}\|^2,$  $\int_0^l|\theta_2(t,l)|^2dt=\int_0^l|\theta_{20}(l-t)|^2dt\leq\|\theta_{20}\|^2.$ Hence $C$ is admissible for $A_0$. Observe that $A_1$ is a bounded linear operator, by \cite{Weiss1989}, $C$ is admissible for $A$.
The transfer function of delay free system is given by
\begin{align*}
  G(s)=\frac{1}{h_1+h_2}\left(
                          \begin{array}{cc}
                            h_2e^{-sl}-h_2e^{-(h_1+h_2+s)l} & h_2e^{-(h_1+h_2+s)l}+h_1 e^{-sl} \\
                            h_2e^{-sl}+h_1e^{-(h_1+h_2+s)l} & -h_1e^{-(h_1+h_2+s)l}+h_1 e^{-sl} \\
                          \end{array}
                        \right)
\end{align*}
Obviously, $G(s)$ is bounded on $Re s>0$. This implies by \cite{Weiss1994} that the delay free system corresponding to (\ref{delay}) is a well-possed, that is, the current state and the output are continuously depended on the initial state and input. By the same procedure as \cite[Theorem 2.1]{Guo2009}, for system (\ref{delay}) with $\tau=0$,  {\it the input belonging $L^2(\tau,\infty)\times L^2(\tau,\infty)$ implies the output belonging $L^2(0,\infty)\times L^2(0,\infty)$}.
This is very important, because the output is considered as an input in the observer design.
Moreover $G(s)\rightarrow 0$ as $s\rightarrow +\infty$. This implies by \cite{Weiss1994} that $A-BKC_\Lambda$ is a generator of $C_0$-semigroup and $B$ is admissible for $A-BKC_\Lambda$, where $K=\left(
                            \begin{array}{cc}
                              k_1 & 0 \\
                              0 & k_2 \\
                            \end{array}
                          \right)
$ and $C_\Lambda$ being some extension of operator $C$. Moreover, $A-BKC_\Lambda$ is system operator with feedback
$\big(u_1(t),u_2(t)\big)^T =-K\big(y_1(t),y_2(t)\big)^T$.

We design the observer and predictor as follows:\\
\textbf{Observer}:  Because of the existence of delay, the state
$\{\big(\theta_1(s,x),\theta_2(s,x)\big), s\in [0,t-\tau], t>\tau\}$ should be estimated from the known observation $\{\big(y_1(s+\tau),y_2(s+\tau)\big),s\in[0,t-\tau],t>\tau\}$. A Luenberger observer is designed by
\begin{align}\label{observe}
    \left\{
      \begin{array}{lll}
        \frac{\partial}{\partial s}\hat{\theta}^+_1(s,x)=-\frac{\partial}{\partial x}\hat{\theta}^+_1(s,x)+h_1(\hat{\theta}^+_2(s,x)-\hat{\theta}^+_1(s,x)), 0<x<l, 0<s<t-\tau, & \hbox{ } \\
        \frac{\partial}{\partial s}\hat{\theta}^+_2(s,x)=-\frac{\partial}{\partial x}\hat{\theta}^+_2(s,x)+h_2(\hat{\theta}^+_1(s,x)-\hat{\theta}^+_2(s,x)), 0<x<l, 0<s<t-\tau, & \hbox{ } \\
        \hat{\theta}^+_1(s,0)=-k_1[\hat{\theta}^+_2(s,1)-y_1(s+\tau)]+u_1(s), s\geq 0,& \hbox{ } \\
        \hat{\theta}^+_2(s,0)=-k_2[\hat{\theta}^+_2(t,1)-y_2(s+\tau)]+u_2(s), s\geq 0,& \hbox{ } \\
        \hat{\theta}^+_1(0,x)=\hat{\theta}_{10}(x),\hat{\theta}^+_2(0,x)=\hat{\theta}_{20}(x),& \hbox{ }
      \end{array}
    \right.
\end{align}
where $\hat{\theta}_{10}(x)$ and $\hat{\theta}_{20}(x)$ are the (arbitrarily assigned)
 initial state of the observer. System (\ref{observe}) can be written as the form $z(s)=(A-BKC_\Lambda)z(s)+B\big[\big(u_1(s),u_2(s)\big)^T-K\big(y_1(s+\tau),y_2(s+\tau)\big)^T\big]$. The
 admissibility of $B$ for $A-BKC_\Lambda$ implies that the current state of system (\ref{observe}) depends continuously on the initial state and the $L^2$ norm of $\big(u_1(\cdot),u_2(\cdot)\big)^T$ and $\big(y_1(\cdot+\tau),y_2(\cdot+\tau)\big)^T$.\\
 \textbf{Predictor}: Predict $\{\big(\theta_1(s,x),\theta_2(s,x)\big),\ s\in(t-\tau,t],t>\tau\}$ by  $\{\big(\hat{\theta}^{-}_1(s,t,x),\hat{\theta}^{-}_2(s,t,x)\big),\ s\in[0,t-\tau],t>\tau\}$.
For this purpose, we solve (\ref{delay}) with estimated initial values $\big(\hat{\theta}^+_1(t-\tau,x),\hat{\theta}^+_2(t-\tau,x)\big)$ obtained from (\ref{observe})
\begin{align}\label{observer}
    \left\{
      \begin{array}{lll}
        \frac{\partial}{\partial s}\hat{\theta}^{-}_1(s,t,x)=-\frac{\partial}{\partial x}\hat{\theta}^{-}_1(s,t,x)+h_1(\hat{\theta}^{-}_2(s,t,x)-\hat{\theta}^{-}_1(s,t,x)), 0<x<1, t-\tau<s<t,& \hbox{ } \\
        \frac{\partial}{\partial s}\hat{\theta}^{-}_2(s,t,x)=-\frac{\partial}{\partial x}\hat{\theta}^{-}_2(s,t,x)+h_2(\hat{\theta}^{-}_1(s,t,x)-\hat{\theta}^{-}_2(s,t,x)), 0<x<1, t-\tau<s<t, & \hbox{ } \\
        \hat{\theta}^{-}_1(s,t,0)=u_1(s), t-\tau\leq s\leq t,& \hbox{ } \\
        \hat{\theta}^{-}_2(s,t,0)=u_2(s), t-\tau\leq s\leq t,& \hbox{ } \\
        \hat{\theta}^{-}_1(t-\tau,t,x)=\hat{\theta}^+_1(t-\tau,x),\hat{\theta}^{-}_2(t-\tau,t,x)=\hat{\theta}^+_2(t-\tau,x),0\leq x\leq l. & \hbox{ }
      \end{array}
    \right.
\end{align}

With the above {\it observer} and {\it predictor}, we design the {\it estimated state feedback control law} by
\begin{align*}
    u_1=\left\{
          \begin{array}{ll}
            k_1 \hat{\theta}^{-}_1(t,t,l), t>\tau,& \hbox{} \\
            0,  0\leq t\leq \tau, & \hbox{}
          \end{array}
        \right.
     u_2=\left\{
          \begin{array}{ll}
            k_2 \hat{\theta}^{-}_2(t,t,l), t>\tau, & \hbox{} \\
            0, 0\leq t\leq \tau. & \hbox{}
          \end{array}
        \right.
\end{align*}
Denote $\varepsilon^+_1(s,x)=\hat{\theta}^+_1(s,x)-\theta_1(s,x), \varepsilon^+_2(s,x)=\hat{\theta}^+_2(s,x)-\theta_2(s,x),\ 0\leq s\leq t-\tau;$ $\varepsilon^{-}_1(s,t,x)=\hat{\theta}^{-}_1(s,t,x)-\theta_1(s,x), \varepsilon^{-}_2(s,t,x)=\hat{\theta}^{-}_2(s,t,x)-\theta_2(s,x), t-\tau\leq s\leq t$.
The closed-loop system is transferred to the following partial differential equations
\begin{align}\label{del}
    \left\{
      \begin{array}{lll}
        \frac{\partial}{\partial t}\theta_1(t,x)=-\frac{\partial}{\partial x}\theta_1(t,x)+h_1(\theta_2(t,x)-\theta_1(t,x)), 0<x<l, t>\tau, & \hbox{ } \\
        \frac{\partial}{\partial t}\theta_2(t,x)=-\frac{\partial}{\partial x}\theta_2(t,x)+h_2(\theta_1(t,x)-\theta_2(t,x)), 0<x<l, t>\tau, & \hbox{ } \\
        \theta_1(t,0)=-k_1\varepsilon^{-}_2(t,t,l)-k_1\theta_2(t,l), t>\tau,& \hbox{ } \\
        \theta_2(t,0)=-k_2\varepsilon^{-}_1(t,t,l)-k_2\theta_1(t,l), t>\tau,& \hbox{ }
      \end{array}
    \right.
\end{align}
\begin{align}\label{efuu}
    \left\{
      \begin{array}{lll}
        \frac{\partial}{\partial s}\hat{\varepsilon}^+_1(s,x)=-\frac{\partial}{\partial x}\hat{\varepsilon}^+_1(s,x)+h_1(\hat{\varepsilon}^+_2(s,x)-\hat{\varepsilon}^+_1(s,x)), 0<x<l, 0<s<t-\tau, t>\tau, & \hbox{ } \\
        \frac{\partial}{\partial s}\hat{\varepsilon}^+_2(s,x)=-\frac{\partial}{\partial x}\hat{\varepsilon}^+_2(s,x)+h_2(\hat{\varepsilon}^+_1(s,x)-\hat{\varepsilon}^+_2(s,x)), 0<x<l, 0<s<t-\tau, t>\tau, & \hbox{ } \\
        \hat{\varepsilon}^+_1(s,0)=-k_1\hat{\varepsilon}^+_2(s,1), 0\leq s\leq t-\tau, t>\tau,& \hbox{ } \\
        \hat{\varepsilon}^+_2(s,0)=-k_2\hat{\varepsilon}^+_2(s,1), 0\leq s\leq t-\tau, t>\tau,& \hbox{ } \\
        \hat{\varepsilon}^+_1(0,x)=\hat{\theta}_{10}(x)-\theta_{10}(x),\hat{\varepsilon}^+_2(0,x)=\hat{\theta}_{20}(x)-\theta_{20}(x), 0\leq x\leq l,& \hbox{ }
      \end{array}
    \right.
\end{align}
\begin{align}\label{observer12}
    \left\{
      \begin{array}{lll}
        \frac{\partial}{\partial s}\varepsilon^{-}_1(s,t,x)=-\frac{\partial}{\partial x}\varepsilon^{-}_1(s,t,x)+h_1(\varepsilon^{-}_2(s,t,x)-\varepsilon^{-}_1(s,t,x)), 0<x<l, t-\tau<s<t, t>\tau, & \hbox{ } \\
        \frac{\partial}{\partial s}\varepsilon^{-}_2(s,t,x)=-\frac{\partial}{\partial x}\varepsilon^{-}_2(s,t,x)+h_2(\varepsilon^{-}_1(s,t,x)-\varepsilon^{-}_2(s,t,x)), 0<x<l, t-\tau<s<t, t>\tau, & \hbox{ } \\
        \varepsilon^{-}_1(s,t,0)=0, t-\tau\leq s\leq t,t>\tau, & \hbox{ } \\
        \varepsilon^{-}_2(s,t,0)=0, t-\tau\leq s\leq t,t>\tau, & \hbox{ } \\
        \varepsilon^{-}_1(t-\tau,t,x)=\varepsilon^+_1(t-\tau,x),\varepsilon^{-}_2(t-\tau,t,x)=\varepsilon^+_2(t-\tau,x), 0\leq x\leq l,t>\tau.& \hbox{ }
      \end{array}
    \right.
\end{align}

Obviously,  the system operator of  (\ref{efuu}) is $A-BKC_\Lambda$. By the same procedure of \cite[Example 2]{Villegas2009}, we can easily prove that there exists positive constants $M$ and $\gamma_0$ such that
\begin{align}\label{ABKC}
    e^{(A-BKC_\Lambda)t}\leq Me^{-\gamma_0 t},\ \forall \ t\geq 0,
\end{align}
provided $k_1^2<\frac{h_1}{h_2}, k_2^2<\frac{h_2}{h_1}$ hold.
This implies that our observer (\ref{observe}) converges and the inverse operator of $A-BKC_\Lambda$ denoted by $(A-BKC_\Lambda)^{-1}$ exists.
Our main results are described as the following theorem.
\begin{theorem}
Assume that $0<k_1<\sqrt{\frac{h_1}{h_2}},\ 0<k_2<\sqrt{\frac{h_2}{h_1}}$ and $t>\tau$. Then,
for $\tau>l$, system (\ref{del}) decays exponentially for any
initial value; for $0<\tau\leq l$ and
$\left(
   \begin{array}{c}
     w_1 \\
     w_2\\
   \end{array}
 \right)=\left(
   \begin{array}{c}
     \hat{\theta}_{10}-\theta_{10} \\
     \hat{\theta}_{20}-\theta_{20}\\
   \end{array}
 \right)
\in D(A-BKC_\Lambda),$ system (\ref{del}) decays exponentially in the sense that
\begin{align*}
  \bigg\|\left(
    \begin{array}{c}
      \theta_1(t,\cdot) \\
      \theta_2(t,\cdot)  \\
    \end{array}
  \right)\bigg\|
\leq   D_\tau e^{-\gamma t} \bigg[\bigg\|\left(
    \begin{array}{c}
      \theta_{10}(\cdot) \\
      \theta_{20}(\cdot)  \\
    \end{array}
  \right)\bigg\|
+C_\tau^0
\bigg\|
(A-BKC_\Lambda)\left(
  \begin{array}{c}
    w_1 \\
    w_2\\
  \end{array}
\right)\bigg\|\bigg],
\end{align*}
where $D_\tau$ and $C_\tau^0$ are positive constants independent on $t$, $\gamma$ is a given constant smaller than $\gamma_0$.
\end{theorem}

\section{The proof of Theorem 1.1}

\textbf{Proof.} System (\ref{observer12}) tells us $\left(
  \begin{array}{c}
    \varepsilon^{-}_1(t,t,\cdot) \\
    \varepsilon^{-}_2(t,t,\cdot) \\
  \end{array}
\right)=e^{A\tau}
\left(
  \begin{array}{c}
    \varepsilon^{+}_1(t-\tau,\cdot) \\
    \varepsilon^{+}_2(t-\tau,\cdot) \\
  \end{array}
\right).$ By simple computation, we obtain
\begin{align}\label{aaaaa}
    \left\{
      \begin{array}{lll}
       \varepsilon^{-}_1(t,t,x)=&\frac{1}{h_1+h_2}\big[
        (h_2+h_1e^{-(h_1+h_2)\tau})\varepsilon^{+}_1(t-\tau,x-\tau)+h_1(1-e^{-(h_1+h_2)\tau})\varepsilon^{+}_2(t-\tau,x-\tau)\big]\\
        \varepsilon^{-}_2(t,t,x)=& \frac{1}{h_1+h_2}\big[
        h_2(1-e^{-(h_1+h_2)\tau})\varepsilon^{+}_1(t-\tau,x-\tau)+(h_1+h_2e^{-(h_1+h_2)\tau})\varepsilon^{+}_2(t-\tau,x-\tau)\big]
      \end{array}
    \right.
\end{align}
provided $x\geq \tau$, and $\varepsilon^{-}_1(t,t,x)=\varepsilon^{-}_1(t,t,x)=0$, provided $x< \tau$. This implies that, for $\tau>l$, $\varepsilon^{-}_1(t,t,l)=\varepsilon^{-}_2(t,t,l)=0.$
Accordingly, system (\ref{del}) is of the form $\dot{z}(t)=(A-BKC_\Lambda)z(t),t>\tau$, thereby exponentially stable for any initial value.

Below we consider the case $0<\tau\leq l$ and
$\left(
   \begin{array}{c}
     w_1 \\
     w_2\\
   \end{array}
 \right)=\left(
   \begin{array}{c}
     \hat{\theta}_{10}-\theta_{10} \\
     \hat{\theta}_{20}-\theta_{20}\\
   \end{array}
 \right)
\in D(A-BKC_\Lambda)$. In such case, (\ref{aaaaa}) indicates $\varepsilon^{-}_1(t,t,l)=\frac{1}{h_1+h_2}\bigg[
        h_2\varepsilon^{+}_1(t-\tau,l-\tau)+h_1\varepsilon^{+}_2(t-\tau,l-\tau)
        +e^{-(h_1+h_2)\tau}\big(h_1\varepsilon^{+}_1(t-\tau,l-\tau)-h_1\varepsilon^{+}_2(t-\tau,l-\tau)\big)\bigg],$
        $\varepsilon^{-}_2(t,t,l)=\frac{1}{h_1+h_2}\bigg[
        h_2\varepsilon^{+}_1(t-\tau,l-\tau)+h_1\varepsilon^{+}_2(t-\tau,l-\tau)
        +e^{-(h_1+h_2)\tau}\big(-h_2\varepsilon^{+}_1(t-\tau,l-\tau)+h_2\varepsilon^{+}_2(t-\tau,l-\tau)\big)\bigg].$
Hence there exist positive constants $p_1,p_2,q_1,q_2$ independent on $t$ such that
\begin{align}\label{ttl}
   \left\{
     \begin{array}{ll}
       |\varepsilon^{-}_1(t,t,l)|\leq p_1|\varepsilon^{+}_1(t-\tau,l-\tau)|+p_2|\varepsilon^{+}_2(t-\tau,l-\tau)|, & \hbox{} \\
      |\varepsilon^{-}_2(t,t,l)|\leq q_1|\varepsilon^{+}_1(t-\tau,l-\tau)|+q_2|\varepsilon^{+}_2(t-\tau,l-\tau)|. & \hbox{}
     \end{array}
   \right.
\end{align}

Denote
$\left(
   \begin{array}{c}
     \psi_1(t,\tau,\cdot) \\
     \psi_2(t,\tau,\cdot) \\
   \end{array}
 \right)=e^{(A-BKC_\Lambda)(t-\tau)}
(A-BKC_\Lambda)\left(
  \begin{array}{c}
    w_1 \\
    w_2\\
  \end{array}
\right).$
Then, system (\ref{efuu}) tells us that
$\left(
   \begin{array}{c}
     \varepsilon^{+}_1(t-\tau,\cdot) \\
     \varepsilon^{+}_2(t-\tau,\cdot) \\
   \end{array}
 \right)=(A-BKC_\Lambda)^{-1}\left(
   \begin{array}{c}
     \psi_1(t,\tau,\cdot) \\
     \psi_2(t,\tau,\cdot) \\
   \end{array}
 \right)$ and it is of the following form
 \begin{align*}
      \varepsilon^{+}_1(t-\tau,x)
      =&\frac{1}{h_1+h_2}\int_0^x\bigg[h_1e^{(h_1+h_2)(\sigma-x)}\big(\psi_2(t,\tau,\sigma)-\psi_1(t,\tau,\sigma)\big)
      \\
      &-h_2\psi_1(t,\tau,\sigma)-h_1\psi_2(t,\tau,\sigma)\bigg]d\sigma +\int_0^l\big(\alpha_1+\alpha_2e^{(h_1+h_2)\sigma}\big)\psi_1(t,\tau,\sigma)d\sigma \\
      & +\int_0^l\big(\alpha_3+\alpha_4e^{(h_1+h_2)\sigma}\big)\psi_2(t,\tau,\sigma)d\sigma,\\
     \varepsilon^{+}_2(t-\tau,x)
      =&\frac{-1}{h_1+h_2}\int_0^x\bigg[h_2e^{(h_1+h_2)(\sigma-x)}\big(\psi_2(t,\tau,\sigma)-\psi_1(t,\tau,\sigma)\big)d\sigma\\
      &+h_2\psi_1(t,\tau,\sigma)+h_1\psi_2(t,\tau,\sigma) \bigg]d\sigma +\int_0^l\big(\beta_1+\beta_2e^{(h_1+h_2)\sigma}\big)\psi_1(t,\tau,\sigma)d\sigma \\
       &+\int_0^l\big(\beta_3+\beta_4e^{(h_1+h_2)\sigma}\big)\psi_2(t,\tau,\sigma)d\sigma,
 \end{align*}
 where $\alpha_j,\beta_j, j=1,2,3,4$ are nonnegative constants independent on $t$.
Accordingly, there exist positive constants $m_1,m_2,n_1,n_2$ independent on $t$ such that
\begin{align}\label{tl}
   \left\{
     \begin{array}{ll}
       |\varepsilon^{+}_1(t-\tau,l-\tau)|\leq m_1\|\psi_1(t,\tau,\cdot)\|+m_2\|\psi_2(t,\tau,\cdot)\|, & \hbox{} \\
       |\varepsilon^{+}_2(t-\tau,l-\tau)|\leq n_1\|\psi_2(t,\tau,\cdot)+
 n_2\|\psi_2(t,\tau,\cdot)\|. & \hbox{}
     \end{array}
   \right.
\end{align}

We combing (\ref{ttl}) and (\ref{tl}) to get
$|\varepsilon^{-}_1(t,t,l)|\leq (p_1m_1+p_2n_1)\|\psi_1(t,\tau,\cdot)\|+(p_1m_2+p_2n_2)\|\psi_2(t,\tau,\cdot)\|$ and
$|\varepsilon^{-}_2(t,t,l)|\leq (q_1m_1+q_2n_1)\|\psi_1(t,\tau,\cdot)\|+(q_1m_2+q_2n_2)\|\psi_2(t,\tau,\cdot)\|.$
Use (\ref{ABKC}) to obtain
\begin{align}\label{eq}
 \nonumber \bigg| \left(
  \begin{array}{c}
    \varepsilon^{-}_1(t,t,1) \\
    \varepsilon^{-}_2(t,t,1) \\
  \end{array}
\right)\bigg|\leq M_1\bigg\|\left(
   \begin{array}{c}
     \psi_1(t,\tau,\cdot) \\
     \psi_2(t,\tau,\cdot) \\
   \end{array}
 \right)\bigg\|\
 =&M_1\bigg\|e^{(A-BKC_\Lambda)(t-\tau)}
(A-BKC_\Lambda)\left(
  \begin{array}{c}
    w_1 \\
    w_2\\
  \end{array}
\right)\bigg\|\\
&\leq  M_1Me^{-\gamma_0(t-\tau)}\bigg\|
(A-BKC_\Lambda)\left(
  \begin{array}{c}
    w_1 \\
    w_2\\
  \end{array}
\right)\bigg\|,
\end{align}
where $M_1$ is a positive constant depended on $p_1,p_2,q_1,q_2,m_1,m_2,n_1,n_2$.

The abstract form of system (\ref{del}) is described as follows
\begin{align*}
  \frac{d}{dt}\left(
    \begin{array}{c}
      \theta_1(t,\cdot) \\
      \theta_2(t,\cdot)  \\
    \end{array}
  \right)=
  (A-BKC_\Lambda)\left(
    \begin{array}{c}
      \theta_1(t,\cdot) \\
      \theta_2(t,\cdot)  \\
    \end{array}
  \right)+B\left(\begin{array}{c}
    \varepsilon^{-}_1(t,t,1) \\
    \varepsilon^{-}_2(t,t,1) \\
  \end{array}
\right).
\end{align*}
Let $0<\gamma<\gamma_0$ and $Y(t)=e^{\gamma t}\left(
    \begin{array}{c}
      \theta_1(t,\cdot) \\
      \theta_2(t,\cdot)  \\
    \end{array}
  \right)$. We obtain
$$\dot{Y}(t)=(\gamma+A-BKC_\Lambda)Y(t)+Be^{\gamma t}\left(\begin{array}{c}
    \varepsilon^{-}_1(t,t,1) \\
    \varepsilon^{-}_2(t,t,1) \\
  \end{array}
\right), t>\tau $$
with $\{e^{(\gamma+A-BKC_\Lambda)s}\}_{s\geq 0}$ being exponentially stable $C_0$-semigroup.

Since $B$ is admissible for $A-BKC_\Lambda$, $B$ is admissible for $(\gamma+A-BKC_\Lambda)$, see \cite{Weiss1994}. The combination of property of admissibility,  exponentially stability of semigroup $\{e^{(\gamma+A-BKC_\Lambda)s}\}_{s\geq 0}$, and (\ref{eq}) implies
\begin{align*}
  \|Y(t)\|=&\bigg\|e^{(\gamma+A+BC)(t-\tau)}Y(\tau)\|+\int_\tau^t e^{(\gamma+A-BKC_\Lambda)(t-s)}Be^{\gamma s} \left(\begin{array}{c}
    \varepsilon^{-}_1(s,s,1) \\
    \varepsilon^{-}_2(s,s,1) \\
  \end{array}
\right)ds\bigg\| \\
\leq & C_\tau \|Y(\tau)\|+\bigg(\int_\tau^t\bigg\|e^{\gamma s}\left(\begin{array}{c}
    \varepsilon^{-}_1(s,s,1) \\
    \varepsilon^{-}_2(s,s,1) \\
  \end{array}
\right)\bigg\|^2ds\bigg)^{\frac{1}{2}}\\
\leq &  C_\tau \|Y(\tau)\|+MM_2\bigg(\int_\tau^t|e^{\gamma s}e^{-\gamma_0(s-\tau)}|ds\bigg)^{\frac{1}{2}}
\bigg\|
(A-BKC_\Lambda)\left(
  \begin{array}{c}
    w_1 \\
    w_2\\
  \end{array}
\right)\bigg\|,
\end{align*}
where $C_\tau$ and $M_2$ are positive constants independent on $t$.
Hence
\begin{align*}
  \bigg\|\left(
    \begin{array}{c}
      \theta_1(t,\cdot) \\
      \theta_2(t,\cdot)  \\
    \end{array}
  \right)\bigg\|
\leq   D_\tau e^{-\gamma t} \bigg[\bigg\|\left(
    \begin{array}{c}
      \theta_{10}(\cdot) \\
      \theta_{20}(\cdot)  \\
    \end{array}
  \right)\bigg\|
+C_\tau^0
\bigg\|
(A-BKC_\Lambda)\left(
  \begin{array}{c}
    w_1 \\
    w_2\\
  \end{array}
\right)\bigg\|\bigg],
\end{align*}
where $D_\tau=C_\tau e^{\gamma\tau}\|e^{A\tau}\|, C_\tau^0=\frac{MM_0}{C_\tau\sqrt{\gamma_0-\gamma}\|e^{A\tau}\|}$. The proof is therefore completed.



\begin{thebibliography}{1}
\bibitem{Guo2009} B.Z. Guo, K.Y. Yang, Dynamic stabilization of an Euler-Bernoulli beam equation with time delay in boundary observation, Automatica, 2009, 45(6): 1468-1475.
\bibitem{Chen2010} W.Y. Lu, J.H. Chen, Observability of the two-stream parallel-flow heat exchanger equation, IMA Journal of Mathematical Control and Information, 2010, 27(1): 91-102.
\bibitem{Sano2016} H. Sano, Exponential stability of heat exchangers with delayed boundary feedback, IFAC-PapersOnLine, 2016, 49(8): 43-47.
\bibitem{Villegas2009} J.A. Villegas, H. Zwart, Y. Le Gorrec, B. Maschke, Exponential stability of a class of boundary control systems, IEEE Transactions on Automatic Control, 2009, 54(1): 142-147.
\bibitem{Weiss1989} G. Weiss, Admissible observation operators for linear semigroups. Israel J. Mathematics. 1989, 65: 17-43.
\bibitem{Weiss1994} G. Weiss, Regular linear systems with feedback. Math. Control Signals Systems. 1994, 7: 23-57.
\end{thebibliography}
\end{document}